\documentclass[12pt,reqno,a4paper]{amsart}

\usepackage[cp1250]{inputenc}
\usepackage{amsthm}
\usepackage{amsmath}
\usepackage{amsfonts}
\usepackage{bbm}
\usepackage{bbold}
\usepackage{geometry}
\usepackage{graphicx}
\usepackage{ifsym}
\usepackage{xcolor}

\newcommand{\RR}{\mathbbm{R}}
\newcommand{\Rp}{\mathbbm{R}_+}

\newcommand{\Mp}{\mathcal{M}^+(\Rp)}

\newcommand{\BC}{\mathop\mathbf{BC}}		
\newcommand{\C}[1]{{\mathbf{C^{#1}}}}		
\renewcommand{\L}[1]{{\mathbf{L^#1}}}					
\newcommand{\W}[2]{\mathbf{W^{#1,#2}}}		

\newcommand{\f}[1]{ \norm{#1}_{\underline{\mathbf{BC}}} }

\newcommand{\Err}{\operatorname{Err}}

\renewcommand{\d}[1]{\mathinner{\mathrm{d}{#1}}} 
\newcommand{\e}{\mathinner{\mathrm{e}}}		
\renewcommand{\epsilon}{\varepsilon}			
\renewcommand{\phi}{\varphi}					
\newcommand{\modulo}[1]{{\left|#1\right|}}		

\newcommand{\naw}[1]{\left( #1 \right)}
\newcommand{\norm}[1]{\| #1 \|}
\newcommand{\abs}[1]{\left| #1 \right|}

\newtheorem*{ex}{Example}
\let\oldex\ex
\renewcommand{\ex}{\oldex\normalfont}

\newtheorem{defi}{Definition}
\newtheorem{thm}{Theorem}
\newtheorem{lemm}{Lemma}
\newtheorem{remark}{Remark}
\let\oldremark\remark
\renewcommand{\remark}{\oldremark\normalfont}
\newtheorem{proposition}{Proposition}
\newenvironment{proofof}[1]
	{\smallskip\noindent{\textbf{Proof~of~#1.}}
	\hspace{1pt}}{\hspace{-5pt}{\nobreak\quad\nobreak\hfill\nobreak
    	$\square$\vspace{2pt}\par}\smallskip\goodbreak}

\title[Finite Range Approximation Method]{Finite Range Method of Approximation for Balance Laws in Measure Spaces}

\author[P. Gwiazda, P. Orli\'nski, A. Ulikowska]{Piotr Gwiazda$^1$, Piotr Orli\'nski$^1$, Agnieszka Ulikowska$^1$}
\address[$^1$]{
Institute of Applied Mathematics and Mechanics,  University of Warsaw, Warsaw 02-097, Poland}
\email{aulikowska@mimuw.edu.pl}

\begin{document}
\maketitle

{
\centerline{$^1$Institute of Applied Mathematics and Mechanics}
\centerline{University of Warsaw}
\centerline{ul. Banacha 2}
\centerline{02-097 Warsaw, Poland}
} 

\begin{abstract}
\noindent
In the following paper we reconsider a numerical scheme which has been recently introduced in \cite{CGU}. The scheme was designed for a wide class of size structured population models with a nonlocal term describing the birth process. We propose a modification of the algorithm, which eliminates the exponential growth in time of the number of particles constituting a numerical solution. Our approach bases on the Finite Range Approximation of the nonlocal term. We provide a convergence theorem, estimates on the convergence speed and results of numerical simulations for several test cases. 
\end{abstract}

\let\thefootnote\relax

\section{Introduction}

The main purpose of this paper is to present a modification of a numerical scheme introduced in \cite{CGU}. {A purpose of the modification is to prevent the number of particles which constitute a numerical solution from the exponential growth. To perform the task we approximate a nonlocal term, which appears on the right hand side of the model equation \eqref{MainEq}, using the Finite Range Approximation method, see Subsection \ref{FRA} for details.
As shown in \cite{CGU}, a convergence of the scheme follows from the stability estimate \cite[Theorem 2.11 (ii)]{CCGU}. Unfortunately, after the approximation procedure the assumptions of \cite[Theorem 2.11]{CCGU} are not fulfilled anymore. From that reason we need to establish a relaxed version of the stability estimate and apply a new strategy in the proof of convergence of the scheme.}

The scheme under consideration follows a current trend which bases on a kinetic approach to population dynamics problems \cite{Ackleh1, Ackleh2, Muntean, siamania, GwiazdaThomasEtAl, GwiazdaMarciniak, piccoli}. Within this approach a population of individuals is divided into groups, which are called cohorts. In other words, a distribution of the population is approximated by a sum of Dirac measures, each one of which represents the average state and the number of individuals within the corresponding cohort.
Such a method of approximation is vastly convenient for numerical studies, especially when it comes to compatibility of a model with an experimental data. Indeed, a result of a measurement of a population is usually a number of individuals which state is within a specified range. A good example of such measurements are demographical studies which provide data about a size of age-cohorts.

A broad group of methods originated from the kinetic theory are particle methods, which are designed to model a behavior of large groups of interacting particles or individuals. Over the last decades they have been successfully applied to solve numerically many problems originated from physics as the Euler equation in fluid mechanics \cite{Euler2, Euler} and the Vlasov equation in plasma physics \cite{BirdLang,CR,GV}. Recently, the particle methods have been used in problems related to crowd dynamics and flow of pedestrians \cite{Muntean,piccoli, PT}, models of a collective motion of large groups of agents \cite{CCR, Dorsogna, KSUB} and population dynamics \cite{CCGU}. For more applications see \cite{Harlow,Issautier, Raviart,Tadmor} and references therein.

In this paper we focus on the population dynamics and the following size structured population model


%


\begin{align}\label{MainEq}
\frac{\partial}{\partial t} \mu + \frac{\partial}{\partial
x}(b(t,\mu)\mu) + c(t,\mu) \mu & = \int_{\Rp}(\eta(t,\mu)
)(y)\d \mu(y),
\end{align}
where $t \in [0,T]$ and $x \geq 0$ denote, respectively, time and the size of an individual. In general, the $x$ variable can describe other physiological states (e.g. length or weight) but for sake of simplicity we stick to the size variable. The measure $\mu$ is a distribution of individuals with respect to $x$.  We assume that an individual changes its size according to the following ODE
\begin{equation}\label{beee}
 \dot{x} = b(t, \mu)(x),
\end{equation}
where $b$ describes a dynamics of the transformation, that is, a speed of the individual's growth. $c(t,\mu)(x)$ simply denotes a death rate, and the integral term describes a birth process. Let us briefly explain the meaning of the right hand side of \eqref{MainEq}. For simplicity we assume for a moment that the $\eta$ function does not depend on time $t$ nor the population state $\mu$. Then, for a fixed $y \geq 0$, $\eta(y)$ describes a distribution (with respect to $x$) of offsprings of an individual of the size $y$. Therefore, the integral in \eqref{MainEq} describes a distribution of all new born individuals at each time moment. In the particular case where all new born individuals have the same size $x^b$ we set
\begin{equation}\label{form_of_eta3}
\eta(y) = \beta(y) \delta_{x = x^b}\,,
\end{equation}
where $\beta(y)$ is related to the probability that an individual of the size $y$ procreates. If \eqref{form_of_eta3} holds, then the integral in \eqref{MainEq} transforms into a boundary condition and, as a consequence, \eqref{MainEq} can be reduced to the following classical renewal equation with the nonlocal boundary condition
\begin{eqnarray}\label{MainEq_b}
\frac{\partial}{\partial t} \mu + \frac{\partial}{\partial
x}(b\mu) + c\mu & = & 0,
\\
\nonumber
b(x^b) D_{\lambda}\mu(x^b) & = & \int_{\Rp}\beta(y)\d \mu(y),
\end{eqnarray}
where $D_{\lambda}\mu(x^b)$ is the Radon-Nikodym derivative of $\mu$ with respect to the Lebesgue measure at $x^b$.

The model \eqref{MainEq} describes a population which comes under processes of birth, death and development. A number of individuals in the population and its total biomass change in time, which clearly indicates a nonconservative character of the problem. We need to underline that the lack of conservativity is the main challenge associated with an application of the particle methods in the population dynamics.
Let us briefly recall that the most common mathematical framework for the kinetic theory is a space of probability measures equipped with a Wasserstein distance. Unfortunately, the {$1$-Wasserstein distance $W_1$} between two Radon measures $\mu$ and $\nu$ such that $\int \d \mu \neq \int \d\nu$ is infinite, which is the reason why natural distances for measures, like the Wasserstein distances, cannot be exploited in case of the nonconservative problems.
{
Indeed, let $\mu,\nu$ be finite Radon measures on $\RR$ such that $\mu(\RR) \neq \nu(\RR)$. Then, according to \cite[Definition 6.1]{villani}
\begin{eqnarray}\label{waser}
W_1(\mu, \nu) 
&=&
 \sup \left(
\int_{\RR} \phi(x)\d (\mu - \nu)(x) \; : \; \mathbf{Lip}(\phi) \leq 1
\right)
\\
\nonumber
&\geq&
 \sup_{a \in \RR} \int_{\RR} a\; \d (\mu - \nu)(x) = a (\mu(\RR) - \nu(\RR)),
\end{eqnarray}
and thus $W_1(\mu,\nu) = +\infty$.}
Therefore, a suitable framework which allows to establish well-posedness of the population dynamics models in the space of measures has to be developed at the first place.
This has been recently achieved by replacing the Wasserstein distance by the {flat metric} (see Section~\ref{mes} for definitions and technical details).

One of the first steps in that field has been made in \cite{GwiazdaThomasEtAl, GwiazdaMarciniak}, where existence, uniqueness and stability of solutions to \eqref{MainEq_b} in the space of finite, nonnegative Radon measures equipped with the flat metric were proved. Within the latter framework the first formal proof of convergence of a corresponding particle method for  \eqref{MainEq_b} has been conducted in \cite{EBT}. The method is called the Escalator Boxcar Train (EBT), and although it was described for the first time in 80's in \cite{deRoos}, the proof of its convergence and the convergence rate \cite{GJMU} is very recent.

Well posedness of a general size-structured population model \eqref{MainEq} in the space of measures was established in \cite{CCGU}, and a numerical scheme based on the particle methods was developed in \cite{CGU}.
In the latter paper an essential assumption is the particular form of the $\eta$ function, namely
\begin{equation}\label{form_of_eta}
\eta(y) = \sum_{p=1}^{r}\beta_p(y) \delta_{x = f_p(y)},
\end{equation}
which means that the size of a child belongs to a set $\{f_p(y)\}$, $p = 1, \dots,
r$, where $y$ is its parent size. For instance, setting $r=1$, $f_1(y) = x^b$, and {$\beta_1(y) = \beta > 0$} corresponds to the special case of \eqref{form_of_eta3}, and leads to the equation \eqref{MainEq_b}.
Another common example is a simple symmetric cell division model, which arises by setting $r = 1$, $f_1(y) = \frac{1}{2}y$, and {$\beta_1(y) = \beta > 0$}. The asymmetric case is obtained by setting $r=2$, $f_1(y) = \sigma y$, $f_2(y) = (1 - \sigma)y$, where $0 < \sigma < 1$, {and $\beta_1(y) = \beta_1 > 0$, $\beta_2(y) = \beta_2 > 0$}.

As it has been already stated above, in the kinetic approach a solution is approximated by a sum of Dirac measures at each discrete time moment. In case of the algorithm developed in \cite{CGU} Dirac deltas represent cohorts, that is groups of individuals of a similar size. Since a population comes under a process of births,  at least one additional Dirac measure is created at each time step of the algorithm. In case of the equation \eqref{MainEq_b} it is exactly one Dirac measure, since all new born individuals have the same size. However, in case of the symmetric cell division model the number of Dirac measures is doubled at each time step, which results in the exponential growth of particles. In order to prevent this phenomenon authors of \cite{CGU} developed a reconstruction procedure, which is in fact an approximation procedure. More precisely, if too many Dirac deltas are created at a particular time step, they are simply approximated by a measure composed of a smaller number of them. Nevertheless, the reconstruction has to be performed once per several time steps, which influences an accuracy of a numerical solution.

Our improvement of the scheme presented in \cite{CGU} bases on a new way of approximation. We postulate to approximate properly the $\eta$ function given by \eqref{form_of_eta} before performing any numerical simulation. Namely, we approximate  $f_p$ functions in \eqref{form_of_eta} by piecewise constant functions $f_p^{\epsilon}$ and, as a consequence, new Dirac measures appear only at some fixed points of the ambient space. Therefore, we run the scheme with a bit inaccurate coefficient, but it turns out that we do not need to perform any approximations nor reconstructions during its execution.  The main problem of such approximation is that the new $\eta$, that is
$$
\tilde \eta(y) = \sum_{p=1}^{r}\beta_p(y) \delta_{x = \tilde f_p(y)},
$$
does not fulfill the assumptions \cite[Assumptions (3.1) - (3.4)]{CGU} providing well posedness of \eqref{MainEq} in the space of measures, which is caused by a fact that the piecewise constant functions are not a subclass of the Lipschitz functions. Fortunately, in the following paper we have overcome this obstacle by developing a relaxed version of the stability estimate, see Subsection \ref{weak} (Remark \ref{f_bar} and Theorem \ref{better_estimate}). To accomplish the task we assumed that $f(x) \leq x$, which is not a restrictive limitation since the sublinearity of $f$ is biologically justified. {Let us mention just the basic examples, i.e. the age of a new born individual is always equal to zero, which is not greater that its parent's age, the sizes of daughter cells are smaller than the size of a mother cell before the mitosis process, a polymer chain is shorter after the division process. For more examples see \cite[Section~3]{CCGU}}.

This paper is organized as follows. In Section \ref{sec_opis} we briefly describe the numerical scheme and the Finite Range Approximation method. Section~\ref{mes} consists of some basic facts about the space of finite, nonnegative Radon measures equipped with the flat metric. For the sake of completeness we also justify the choice of the latter space instead of the Banach space $(\W{1}{\infty})^*$. In Section \ref{conv} we provide convergence results. In Section \ref{num} we show results of numerical simulations for several test cases.

\section{Splitting-Particle Method with Finite Range Approximation}\label{sec_opis}
\subsection{Finite Range Approximation}\label{FRA}
In the following subsection we show how to construct the Finite Range Approximation of a Lipschitz continuous function. The Lipschitz continuity of $f_p$ in \eqref{form_of_eta} is the assumption required for well-posedness of \eqref{MainEq}, see \cite[Assumption (3.4)]{CGU}. Before we proceed let us take a closer look at the following example.
\begin{ex}
Setting $r=1$, $f_1(y) = \frac{1}{2}y$, and {$\beta_1(y) = \beta > 0$} in \eqref{form_of_eta} yields the symmetric cell division model. Application of the particle-based scheme developed in \cite{CGU} to that model results in the exponential growth of Dirac measures approximating a solution. It is a consequence of the fact that a child's size is exactly a half of its parent size. If we substitute the function $f_1$ by a suitable piecewise constant approximation, then a set of all possible sizes of the children becomes finite.
\begin{figure}[h]
\centering
\includegraphics[width=260px]{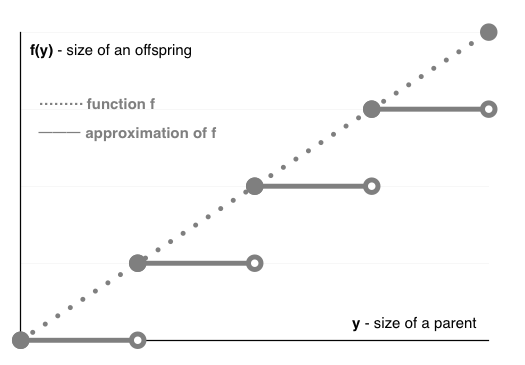}
 \caption{The figure shows the function $f_1$ and its Finite Range Approximation.}
 \label{figure:example}
\end{figure}
\end{ex}
\begin{defi} Let $f : \RR \to \RR$ be a Lipschitz function. We said that $f^{\epsilon}:\RR \to \RR$ is the Finite Range Approximation of $f$ on the interval $[0, M)$, if
$$
\# {\rm{Im}(f^{\epsilon})} < \infty,
\quad\mathrm{and}\quad
\norm{(f - f^{\epsilon})|_{[0,M)}}_{\L\infty} < \epsilon,
$$
where $\epsilon, M > 0$ are arbitrary constants, $\rm{Im}(f^{\epsilon})$ is the image of the function $f^{\epsilon}$, and $\# A$ denotes the number of elements in a set $A$.
\end{defi}
\noindent
The next lemma shows how to construct the Finite Range Approximation of a Lipschitz continuous function.
\begin{lemm}\label{ccc}
Let $\epsilon, M > 0$, and $f: \RR \to \RR$ be a Lipschitz continuous function such that $f(x) \leq x$, for all $0 \leq x < M$. Then, there exists the Finite Range Approximation of $f$, that is a function $f^{\epsilon}$, such that $f^{\epsilon}(x) \leq x$.
\end{lemm}

\begin{proofof}{Lemma \ref{ccc}}
Since $f(x) \leq x$, it holds that $f|_{[0,M)} < M$. Let
\[ J = \left\{
  \begin{array}{l l}
    \left\lfloor \frac{M}{\epsilon} \right\rfloor, & \quad \text{if $\frac{M}{\epsilon} \in \mathbbm{N},$}
\\[2mm]
    \left\lfloor \frac{M}{\epsilon} \right\rfloor +1, & \quad \text{if $\frac{M}{\epsilon} \notin \mathbbm{N}$}.
  \end{array} \right.\]
Define $A_{j} = f^{-1}(\;[(j-1)\epsilon, j\epsilon)\;)$, for $j=1, \dots, J$. It follows directly from the construction that $\cup_{j=1}^{J}A_{j}= [0,M)$. Since $f^{\epsilon}$ is supposed to be defined for all $x \in \RR$, we redefine the sets
$A_1 := (-\infty, 0) \cup A_1$, and $A_{J} := A_{J} \cup [M,+\infty)$. The approximation $f^{\epsilon}$ is thus given by the formula
\begin{equation}\label{ef}
f^{\epsilon}(x) = \sum_{j=1}^{J} a_j \chi_{A_{j}}(x),
\end{equation}
where $a_j = (j-1)\epsilon$ and $\chi_{A_{j}}$ is the characteristic function of the set $A_{j}$. It follows directly form the construction of $f^{\epsilon}$ that it is the Finite Range Approximation of $f$ and $f^{\epsilon}(x) \leq x$, for all $0 \leq x < M$.
\end{proofof}
\begin{remark}\label{support}
For our purposes, it is sufficient to consider a finite interval $[0,M)$. Note that the equation
\begin{eqnarray}\label{h_eq}
\frac{\partial}{\partial t} \mu + \frac{\partial}{\partial
x}(b(t,\mu)\mu) + c(t,\mu) \mu &=& 0,
\end{eqnarray}
admits a finite propagation speed property. Assume that the support of the initial data $\mu_o$ is contained in the interval $[0,M_o)$, for some $M_o > 0$. Then, the support of a solution $\mu(t)$ at time $t \in [0,T]$ is a subset of the interval $[0,M)$, for some $M > 0$. The constant $M$ depends on $M_o$, a suitable norm of $b$, and the length of the time interval $[0,T]$. 
Substituting the right hand side of \eqref{h_eq} by $\int \eta \d \mu$, where $\eta$ is given by \eqref{form_of_eta} and assuming that $f_p(x) \leq x$, for $p=1,\dots,r$, do not influence the support of $\mu(t)$, for any $t \in [0, T]$.
\end{remark}

\subsection{Description of the algorithm}\label{description}

Fix $T>0$ and $N \in \mathbbm{N}$. Define the length of a time step $\Delta t = T/N$ and a set of discrete time points $t_k = k \Delta t$, where $k = 0, \dots, N$.
%
%
Assume that a numerical solution $\mu_k$ at time $t_k$ is a sum of $N_k$ Dirac deltas, that is,
$$
\mu_{k} = \sum_{i=1}^{N_k}m^i_k \delta_{x^i_k},
$$
where $m^i_k$ is a mass of the $i$-th Dirac delta and $x^i_k$ denotes its location. In particular, we assume that the initial data is a sum of Dirac measures. If it is not the case, the initial data can be approximated by such a sum, see Lemma~\ref{init} for details.
Let $M > 0$ be a constant, such that the support of a solution to \eqref{MainEq} is contained in the interval $[0,M)$, for all $t \in [0,T]$. Such a constant exists according to Remark~\ref{support}.
Since the algorithm bases on the splitting technique, it is divided into two steps described below.
\\[1mm]
\noindent
\textbf{Step 1.} The first step is to calculate the new locations of Dirac measures. It is obtained by solving the following ODEs system
\begin{equation}\label{subproblem1}
\frac{d}{dt} x^i(t) = b_k(x^i(t)),\quad
x^i(t_k) = x^i_k,\quad i = 1, \dots, N_k,
\quad
\mathrm{where}
\quad
b_k = b(t_k, \mu_k),
\end{equation}
on a time interval $[t_k, t_{k+1}]$. A result of this step is a measure
$$\bar \mu_{k} = \sum_{i=1}^{N_k}m^i_k\delta_{x^i_{k+1}},
\quad 
\mathrm{where}
\quad
x^i_{k+1} = x^i(t_{k+1}).
$$
\textbf{Step 2.} The second step is to determine locations of the new Dirac measures, which correspond to the newly born individuals. We also need to recalculate masses of all Dirac measures. In order to determine the locations of the new particles we use the Finite Range Approximation of functions $f_p$ in \eqref{form_of_eta} provided by Lemma \ref{ccc}.
According to \eqref{ef}, new particles appear only at points $\{a_j\}_{j=1}^J$, and these locations are fixed in time. As a consequence, we obtain a set of $N_{k+1} := N_k + J$ Dirac measures. In order to recalculate the masses we need to solve the following system of ODEs
\begin{equation}\label{subproblem2}
\left\{
\begin{array}{rcl}
\displaystyle\frac{d}{dt} m^i(t) &=& c_k(x^i_{k+1})m^i(t),
\\[2mm]
 m^i(t_k) &=& m^i_{k}, \quad 1 \leq  i  \leq N_k,
\\[3mm]
\displaystyle\frac{d}{dt} m^i(t) &=& c_k(x^i_{k+1})m^i(t) + \sum_{j=1}^{N_{k+1}}\sum_{p=1}^{r} m^j(t) \alpha_{pk}(x_{k+1}^j, x_{k+1}^i),
\\[2mm]
m^i(t_k) &=& 0,
\quad
N_k + 1 \leq i \leq N_{k+1},
\end{array}
\right.
\end{equation}
where $c_k = c(t_k,\bar \mu_k)$, $\beta_{pk} = \beta_p(t_k,\bar \mu_k)$
and 
$$
\quad\alpha_{pk}{(x_{k+1}^j, x_{k+1}^i)} = \left\{ 
   \begin{array}{l l}
     \beta_{pk}(x_{k+1}^j), & \quad \text{if}\;\;  f^{\epsilon}_p(x_{k+1}^j) = x_{k+1}^i,\\[2mm]
     0, & \quad \text{otherwise.} 
   \end{array} \right.
$$
%
%
After indexes reassignment (e.g., arranging all Dirac measures in the ascending order with respect to the location) we obtain a measure $\mu_{k+1} = \sum_{i=1}^{N_{k+1}}m^i_{k+1}\delta_{x^i_{k+1}}$, which is the output of the algorithm at time $t_{k+1}$. The essential feature of this algorithm is that after $k$ time steps the number of Dirac deltas approximating the solution is equal to $N_o + Jk$, where $N_o$ is the initial number of Dirac measures.
\\[3mm]

\section{The Space of Measures $(\Mp, \rho_F)$}\label{mes}


Henceforth, $\Mp$ denotes the space of nonnegative Radon measures with bounded total variation on $\Rp := \{x \in \mathbbm R \;\colon x \geq 0\}$.  
We equip $\Mp$ with the flat metric
\begin{equation}
  \label{distance}
  \rho_F(\mu_1, \mu_2)
  =
  \sup
  \left\{
    \int_{\Rp} \phi \, \d(\mu_1-\mu_2)
    \;\;\colon\;\;
    \phi \in \C{1}(\Rp)
    \;\;\mbox{and}\;\;
    \norm{\phi}_{\W{1}{\infty}} \leq 1
  \right\},
\end{equation}
where $\norm{\phi}_{\W{1}{\infty}} = \max \left\{
\norm{\phi}_{\L\infty}, \norm{\partial_x \phi}_{\L\infty}
\right\}$. 
\noindent
The condition $\C1(\Rp)$ in \eqref{distance} can be replaced by  $\W{1}{\infty}(\Rp)$ through a standard mollifying sequence argument applied to the test function $\phi$, as its derivative is not involved in the value of the integral, which implies that $\rho_F$ is the metric dual to the $\norm{\cdot}_{({\W{1}{\infty}})^{*}}$ distance. Note that in this paper, the space $\Mp$ is equipped with the metric $\rho_F$ and this shall remain until said differently.
The space $(\Mp,\rho_F)$ is complete and separable.
Using the standard mollification procedure, $\rho_F$ can be equivalently rewritten as
\begin{equation*}
  \label{distance2}
  \rho_F(\mu_1, \mu_2)
  =
  \sup
  \left\{
    \int_{\Rp} \phi \, \d(\mu_1-\mu_2)
    \;\;\colon\;\;
    \phi \in \W{1}{\infty}(\Rp)
    \;\;\mbox{and}\;\;
    \norm{\phi}_{\W{1}{\infty}} \leq 1
  \right\},
\end{equation*}
and therefore, for all $\mu_1,\mu_2 \in \Mp$ it holds that $\rho_F(\mu_1, \mu_2) = \norm{\mu_1 - \mu_2}_{(\W{1}{\infty})^*}$. This equality gives a rise to a question about the possibility of setting the model in the Banach space $\naw{\W{1}{\infty},\norm{\cdot}_{(\W{1}{\infty})^*}}$.
It is a natural question, considering the previous papers in which the problem of continuity of solutions with respect to time were addressed. In \cite{webb} authors established the continuity of solutions to the age structured population model \eqref{MainEq_b} in the $\L{1}$ topology. Then, it was proved in \cite{diekmann} that solutions to balance laws in the space of measures are continuous in the weak-* topology of the Radon measures space. 
Unfortunately, it turns out that the answer to our initial question is negative, because shift operators are not continuous in $\naw{\W{1}{\infty},\norm{\cdot}_{(\W{1}{\infty})^*}}$ in contrary to the $\L1$ topology. This is the obstacle one cannot overcome, since the continuity of the shift operators is essential for obtaining the continuity of solutions to the transport equation.


\begin{lemm}\label{Lema1}
One-parameter semigroup $\{T\}_{t\geq 0}$ of the shift operators is not a strongly continuous semigroup on the space $(\W{1}{\infty}(\Rp))^*$. Moreover, for all $t > 0$ it holds that 
\begin{equation*}
\norm{T_t - I}_{L\left((\W{1}{\infty})^* \right)} \geq 1.
\end{equation*}
{where $L((\W{1}{\infty})^*)$ is the operator norm.}
\end{lemm}
\noindent
The proof of Lemma \ref{Lema1} can be found in Appendix.

\section{Convergence of the Algorithm}\label{conv}

\subsection{Theoretical results concerning well-posedness of \eqref{MainEq}}
In this subsection we recall theoretical results from \cite{CCGU} concerning well posedness of \eqref{MainEq}. Assume that
\begin{align}
\nonumber
\label{eqAssumptions}
& b,c, \beta_p \in {\BC}^{\mathbf{\alpha,1}}
\left(
[0,T] \times {\mathcal M}^{+}(\Rp); \;\W{1}{\infty}(\Rp)
\right),\;\mathrm{for}\;\;p = 1, \dots, r,
\\[1mm]
& f_p \in \mathbf{Lip}(\Rp; \Rp), \;\; f_p(x)\leq x,\;\mathrm{for}\;\;p = 1, \dots, r,
\\[1mm]
\nonumber
& b(t, \mu)(0) \geq 0,\;\;\mathrm{for}\;\; (t,\mu) \in [0,T] \times \Mp.
\end{align}
Here, ${\BC}^{\mathbf{\alpha,1}}([0,T]\times{\mathcal
M}^+(\Rp); \W{1}{\infty}(\Rp))$ is the space of
$\W{1}{\infty}(\Rp)$ valued functions which are
bounded in the $\norm{\cdot}_{\W{1}{\infty}}$ norm, H\"older
continuous with exponent $0<\alpha\leq 1$ with respect to time and
Lipschitz continuous in the flat metric $\rho_F$ with respect to the measure
variable. This space is equipped with the
$\norm{\cdot}_{{\BC}^{\mathbf{\alpha,1}}}$ norm defined by
\begin{equation*}
\norm{f}_{\mathbf{\BC^{\alpha,1}}}
= \norm{f}_{\mathbf{BC}} +
\mathbf{Lip}\left(f(t,\cdot)\right) +
\mathrm{H}_\alpha\left(f(\cdot,\mu)\right),
\end{equation*}
where
$$
\norm{f}_{\mathbf{BC}} = \sup_{t\in[0,T], \mu\in{\mathcal{M}^+(\Rp)}} \norm{f(t, \mu)}_{\W{1}{\infty}},
$$
$\mathbf{Lip}(f)$ is the Lipschitz constant of a function $f$ and
\begin{displaymath}
\mathrm{H}_\alpha(f(\cdot,\mu)) :=
\sup_{s_1,s_2\in[0,T]}
\frac{\norm{f(s_1,\mu) - f(s_2,\mu)}_{\W{1}{\infty}}}{\modulo{s_1 - s_2}^{\alpha}} .
\end{displaymath}
To simplify the notation, we define
$$
\norm{x}_{X^n} = \norm{(x_1, \dots, x_n)}_{X^n} := \sum_{i=1}^{n} \norm{x_i}_X,
\quad\mathrm{where}\;\; x = (x_1, \dots, x_n) \in X^n.
$$
Since in \cite{CGU} the specific form \eqref{form_of_eta} of the $\eta$ function has been assumed, we rewrite the original well-posedness theorem \cite[Theorem 2.11]{CCGU} in the terms of $\beta_p$ and $f_p$ instead of $\eta$, see Theorem \ref{thm:Main} below. Regularity of  $\beta_p$ and $f_p$ imposed in
\eqref{eqAssumptions} guarantees
that $\eta$ defined by \eqref{form_of_eta} fulfills the
assumptions of \cite[Theorem 2.11]{CCGU} and thus, \eqref{MainEq} is
well posed. 

\begin{thm}\label{thm:Main}
Let \eqref{eqAssumptions} hold and $\mu_o \in \Mp$.
Then, there exists a unique solution
$$
  \mu \in (\BC\cap\;
  \mathbf{Lip})\Big([0,T] ;({\mathcal M}^+(\Rp),\rho_F)\Big)$$ to~\eqref{MainEq}.
  Moreover, the following properties are satisfied:
  \begin{enumerate}
  \item For all $0 \leq t_1\leq t_2 \leq T$ there exist constants
    $K_1$ and $K_2$, such that
    \begin{equation*}
      \rho_F\left(\mu(t_1),\mu(t_2)\right)
      \leq
      K_1 \e^{K_2({t_2-t_1})} \mu_o(\Rp)({t_2 - t_1}).
    \end{equation*}
  \item Let $\mu_1(0), \mu_2(0) \in {\mathcal M}^{+}(\Rp)$ and
    $b_i$, $c_i$, $\beta_i = (\beta^i_1, \dots, \beta^i_r)$, $f_i = (f^i_1, \dots, f^i_r)$ satisfy
    assumptions \eqref{eqAssumptions} for $i = 1, 2$, $p=1,\dots, r$. Let $\mu_i$
    solve~\eqref{MainEq} with initial datum $\mu_i(0)$ and
    coefficients $(b_i,c_i,\beta_i, f_i)$.  Then, there exist
    constants $C_1$, $C_2$ and $C_3$ such that
    \begin{eqnarray}\label{xxx}
&&    \rho_F\left(\mu_1(t),\mu_2(t)\right)
      \leq \e^{C_1t} \rho_F\left(\mu_1(0),\mu_2(0)\right)
\\[2mm]
&&
\nonumber
\quad\quad\quad\quad+\;
      C_2t \e^{C_3t}\left(
      \norm{
      (b_1-b_2,c_1-c_2,\beta_1-\beta_2)}_{\mathbf{BC}}
  + \norm{f_1- f_2}_{\L\infty}
\right),
    \end{eqnarray}
    for all $t\in [0,T]$.
  \end{enumerate}
\end{thm}
\noindent
All constants in Theorem \ref{thm:Main} depend on suitable norms of the model coefficients. Unfortunately, the constants $C_2, C_3$ in the second claim depend on the Lipschitz constants of $f_1^i, \cdots, f_r^i$, for $i=1,2$. This is the main obstacle we need to overcome, since the Finite Range Approximation of a Lipschitz function is not a Lipschitz function (in fact it is not even continuous). In the next subsection we show how to deal with this problem.

\subsection{Theoretical results concerning well-posedness of \eqref{MainEq} with a relaxed version of the stability estimate}\label{weak}

The first problem with the Finite Range Approximation is that it produces non-continuous functions, which implies that assumptions \eqref{eqAssumptions} are not fulfilled and, as a consequence, there are no results concerning well-posedness of \eqref{MainEq}. At the beginning of this subsection we show how to substitute these noncontinuous functions by suitable Lipschitz continuous functions in such a way that the result of solving \eqref{subproblem1}--\eqref{subproblem2} is not affected. This is a subject of the following remark.

\begin{remark}\label{f_bar}
Define a set $D$ consisting of locations of all Dirac measures at each discrete time moment
$$
D := \{x^i_k,\; i=1,\dots, N_k\}_{k=0}^{N},$$
where $x_k^i$ is the location of the $i$-th Dirac delta at time $t_k = k\Delta t$. This is clearly a finite set, $D \subset [0, \lambda M)$, where $\lambda \geq 1$ and $M$ is the constant from Remark~\ref{support}. 
Let $f_p^{\epsilon}$ be the Finite Range Approximation of $f_p$, for $p=1,\dots,r$. Without loss of generality we assume that $f, f_p^{\epsilon} : \Rp \to \Rp$. For a fixed $p$, $f_p^{\epsilon}$ is a piecewise constant function, which can be rewritten in the following form
$$
f_p^{\epsilon} (x) = \sum_{i=1}^{J_p}a_{i} \chi_{D_{i}^p}(x),
$$
where $a_{i} \geq 0$, $J_p \geq J$, $D_i^p = [d_i^{\;p}, d^{\;p}_{i+1})$, $\cup_{i=1}^{J_p} D_i= \Rp$, and $\{d_i^{\;p}\}_{i=1}^{J_p}$ is a finite, strictly increasing sequence with $d_{J_p+1} = +\infty$. Define
$$
d_i^{\; p, \max} = \max\{x \; : \; x \in D \cap D_i \},\quad \mathrm{for}\;\; i=1,\dots, J_p-1,
$$
and a piecewise linear function $\bar f_p^{\epsilon}$
\begin{equation}\label{F_BARRR}
\bar{f}_p^{\epsilon}(x) = 
\left\{
\begin{array}{lcl}
a_i, 
&\mathrm{for}&
x \in [d_i^{\; p},d_i^{\;p, \max}),\; i < J_p,
\\[1mm]
\displaystyle
\frac{a_{i+1}-a_i}{d^{\;p }_{i+1} - d_{i}^{\; p, max}} x + \frac{a_i d^{\;p}_{i+1} - d_i^{\;p, \max}a_{i+1}}{d^{\;p}_{i+1} - d_i^{\; p,\max}},
&\mathrm{for}&
x \in [d_i^{\; p,\max}, d^{\; p}_{i+1}),\; i < J_p,
\\[1mm]
a_{J_p}, &\mathrm{for}& x \in D_{J_p}.
\end{array}
\right.
\end{equation}
It follows directly form the construction that for all $x \in D$ it holds that
$$\bar f_p^{\epsilon}(x) = f_p^{\epsilon}(x),\quad\mathrm{and}\quad
\norm{ f^{\epsilon}_p - \bar f^{\epsilon}_p}_{\L\infty}\leq \epsilon.$$
\end{remark}
\noindent
Henceforth, we use the modification described above instead of the corresponding Finite Range Approximation. 
Note that the functions $\bar f^{\epsilon}_p$ are Lipschitz continuous, but their Lipschitz constants may increase due to the increase of a number of Dirac measures approximating a solution. Therefore, there is still some work to be done. Namely, we need to obtain an estimate analogous to \eqref{xxx}, but with all constants independent on the Lipschitz constants of the functions  $\bar f^{\epsilon}_p$. This is the subject of the following theorem.

\begin{thm}\label{better_estimate}
Let $\mu_1(0), \mu_2(0) \in {\mathcal M}^{+}(\Rp)$ and
    $b_i$, $c_i$, $\beta_i = (\beta^i_1, \dots, \beta^i_r)$, $f_i = (f^i_1, \dots, f^i_r)$ satisfy
    assumptions \eqref{eqAssumptions} for $i = 1, 2$, $p=1,\dots, r$. Let $\mu_i$
    solve~\eqref{MainEq} with initial datum $\mu_i(0)$ and
    coefficients $(b_i,c_i,\beta_i, f_i)$.  Then, there exists a
    constant $C$, which depends only on suitable norms of $(b_1,c_1,\beta_1, f_1)$, such that
\begin{eqnarray}\label{better_estimate_eq}
&&    \rho_F\left(\mu_1(t),\mu_2(t)\right)
      \leq \e^{ Ct} \rho_F\left(\mu_1(0),\mu_2(0)\right)
\\[2mm]
&&
\nonumber
\quad\quad+\;
       Ct \e^{ Ct}\left(
            \norm{
      (b_1-b_2,c_1-c_2,\beta_1-\beta_2)}_{\underline{\mathbf{BC}}}
+ \norm{f_1 - f_2}_{\L\infty}
\right),
    \end{eqnarray}
where $\norm{f}_{\underline{\mathbf{BC}}} = \sup_{(s,\mu) \in [0,t] \times \Mp}
      \norm{
      f(s,\mu)}_{\L\infty}$.
\end{thm}
\noindent
The power of the estimate \eqref{better_estimate_eq} is that the constant $C$ depends only on coefficients $(b_1, c_1, \beta_1, f_1)$. Therefore, we can plug into  \eqref{better_estimate_eq} the original functions $f_p$ and their approximations $\bar f^{\epsilon}_p$ described in Remark~\ref{f_bar}. Note that the only place where the functions $\bar f^{\epsilon}_p$  appear in \eqref{better_estimate_eq} is the term $\norm{f_1 - f_2}_{\L\infty}$, which can be estimated by $\epsilon$, according to Lemma~\ref{ccc}.
The proof of Theorem \ref{better_estimate} bases on formula \eqref{PropBressan_formula} (see \cite[Theorem 2.9]{bressan_00}), which allows to consider equations locally in time. 
Before we proceed, we introduce some preliminary notions.
\begin{defi}\label{semiflow} Let $(X,\rho)$ be a metric space. A map $S : [0,T] \times X \to X$ is called a Lipschitz semiflow, if
\begin{enumerate}
\item
$S(0,x) = 0$ for all $x \in X$,
\item
$S(t+s,x) = S(t, S(s,x))$ for all $t,s, t+s \in [0,T]$ and $x \in X$.
\item
$\rho(S(t,x), S(s,y)) \leq L_1\rho(x,y) + L_2\modulo{t-s}$.
\end{enumerate}

\end{defi}
\begin{proposition}\label{PropBressan}
Let $S: [0,T] \times X \rightarrow X$ be a Lipschitz semiflow. For every Lipschitz continuous map $\nu:[0,T]\rightarrow X$ the following estimate holds,
\begin{equation}\label{PropBressan_formula}
\rho\big(\nu_t,S(t; 0)\mu_o\big) \leq L_1 \int_{[0,t]} \liminf_{h \downarrow 0}\frac{\rho\big(\nu_{\tau+h}, S(h; \tau)\nu_{\tau}\big)}{h} d{\tau}, 
\end{equation}
where $\rho$ is a corresponding metric.
\end{proposition}
\noindent

\begin{proofof}{Theorem \ref{better_estimate}} Since $C_1$ in \eqref{xxx} depends only on $(b_1, c_1, \beta_1, f_1)$, we may assume that $\mu_1(t)$ and $\mu_2(t)$ are solutions to \eqref{MainEq} with the same initial data $\mu_o$. As the first step, we define the time dependent functions
$$
b_j (t,x) = b_i(t,\mu_j(t))(x),\;\; c_j (t,x) = c_j(t,\mu_j(t))(x),\;\; \beta_{j,p} (t,x) = \beta_p^j(t,\mu_j(t))(x),
$$
for $j=1,2$, $p=1,\dots, r$. Fix $n \in \mathbbm{N}$, define $\Delta t = T/2^n$, $t^i_n = i \Delta t$ for
  $i=0,1,\dots,2^n$, and approximate $b_j$, $c_j$ and $\beta_{j,p}$ as follows:
  \begin{eqnarray*}
    b_n^j (t,x)
    & = &
    \sum_{i=0}^{2^n-1} b_j(t^i_n,x) \, \chi_{\strut [t^{i}_n, t^{i+1}_n)}(t),
    \\
    c_n^j (t,x)
    & = &
    \sum_{i=0}^{2^n-1} c_j(t^i_n,x) \, \chi_{\strut [t^{i}_n, t^{i+1}_n)}(t),
    \\
    \beta^{j,p}_n (t,x)
    & = &
    \sum_{i=0}^{2^n-1} \beta_{j,p}(t^i_n,x) \, \chi_{\strut [t^{i}_n,
      t^{i+1}_n)}(t) \,.
  \end{eqnarray*}
  Note, that on each interval $[t_n^i,t_{n}^{i+1})$ functions defined above do not depend on $t$. Therefore, according to \cite[Theorem 2.8]{CCGU}, solving \eqref{MainEq} with coefficients $b_n^j$, $c_n^j$, and $\beta_n^{j,p}$ on the time interval $[t_n^i,t_{n}^{i+1})$ yields a Lipschitz semigroup. Call $S^{j,i,n}$ the corresponding semigroup
  and define the map $F^{j,n}
 \colon [0,T] \times \Mp \to
  \Mp$ by
  \begin{eqnarray}
    \label{Poligonals:Aut}
    \!\!\!
    F^{j,n}_{t} \mu
    & = &
    \left\{
      \begin{array}{l@{\qquad}l@{}}
        S^{j,0,n}_{t} \mu,
        & \mbox{if }
        t \in [t_n^0, t_n^1),
        \\[10pt]
        \left(
          S^{j,i,n}_{t-t^{i}_n} \circ
          \left(\bigcirc_{q = 0}^{i-1}  S^{j,q,n}_{T/2^n}\right)
              \right)
        \mu,
        & \mbox{if }
        t \in [t_n^{i}, t_n^{i+1}).
      \end{array}
    \right.
  \end{eqnarray}
 It follows from the construction that $F^{j,n}$ is a Lipschitz semiflow. Without loss of generality we assume that $t = t_n^k$. Substituting $F^{1,n}$ and $F^{2,n}$ to \eqref{PropBressan_formula} yields

 \begin{eqnarray*}
 \rho_F\big(F_t^{1,n} \mu_o,F_t^{2,n}\mu_o\big) \leq \mathbf{Lip}(F^{1,n}) \sum_{i=0}^{k-1}\int_{[t_n^i,t_n^{i+1})} \liminf_{h \downarrow 0}\frac{\rho_F\big(S^{1,i,n}_hF^{1,n}_{\tau}\mu_o, S^{2,i,n}_hF^{1,n}_{\tau}\mu_o\big)}{h} d{\tau}.
 \end{eqnarray*}
According to estimates \cite[proof of Theorem 2.8]{CCGU} for the linear autonomous problem, it holds that
\begin{align*}
& \int_{[t_n^i, t_n^{i+1})} \!\!\!\!  \liminf_{h \downarrow 0}  \frac{\rho_F\big(S^{1,i,n}_hF^{1,n}_{\tau}\mu_o, S^{2,i,n}_hF^{1,n}_{\tau}\mu_o\big)}{h} d{\tau}
\\
&\leq
\int_{[t_n^i, t_n^{i+1})} (F^{1,n}_{\tau}\mu_o)(\Rp) d \tau
\Big(\norm{b^1_n(t_n^{i}, \cdot) - b^2_n(t_n^{i}, \cdot) }_{\L\infty} + \norm{c^1_n(t_n^{i}, \cdot) - c^2_n(t_n^{i}, \cdot) }_{\L\infty} 
\\
&\quad\quad+ 
\sum_{p=1}^{r}\norm{\beta^{1,p}_n(t_n^{i}, \cdot) - \beta^{2,p}_n(t_n^{i}, \cdot) }_{\L\infty} + \sum_{p=1}^{r}\norm{f^1_p - f^2_p }_{\L\infty}\Big)
\\
& \leq 
\Big(
 \norm{
      (b_1-b_2,c_1-c_2,\beta_1-\beta_2)}_{\underline{\mathbf{BC}}}
      + \norm{f_1 - f_2}_{\L\infty}
      \Big) \int_{[t_n^i, t_n^{i+1})} (F^{1,n}_{\tau}\mu_o)(\Rp) d \tau.
\end{align*}
Summing over $i=0, \dots, k-1$ yields
\begin{align*}
\rho_F\big(F_t^{1,n} \mu_o,F_t^{2,n}\mu_o\big) & \leq \mathbf{Lip}(F^{1,n}) 
\Big(
 \norm{
      (b_1-b_2,c_1-c_2,\beta_1-\beta_2)}_{\underline{\mathbf{BC}}}
      + \norm{f_1 - f_2}_{\L\infty}
      \Big) 
      \\
      &\quad \cdot \int_{[0, t]} (F^{1,n}_{\tau}\mu_o)(\Rp) d \tau.
\end{align*}
According to the estimate from \cite[proof of Theorem 2.10, claim iii)]{CCGU} it holds that 
$
(F^{1,n}_{\tau}\mu_o)(\Rp) \leq \e^{\bar C\tau}\mu_o(\Rp),
$
where $\bar C$ depends only on $(b_1, c_1, \beta_1, f_1)$. Therefore, there exists a constant $ C$, which depends on the latter set of coefficients and $\mu_o$, such that
\begin{equation}\label{eq123}
 \rho_F\big(F_t^{1,n} \mu_o,F_t^{2,n}\mu_o\big) \leq Ct \e^{Ct}
 \Big(
 \norm{
      (b_1-b_2,c_1-c_2,\beta_1-\beta_2)}_{\underline{\mathbf{BC}}}
      + \norm{f_1 - f_2}_{\L\infty}
      \Big).
\end{equation}
According to \cite{CCGU} (see the proof of Theorem 2.10), the map $F^{j,n}_t\mu_o$ converges uniformly with respect to time to $\mu_j(t)$, as $n\to +\infty$. Therefore, passing to the limit in \eqref{eq123} ends the proof.

\end{proofof}

%

\subsection{Error Estimates in $\rho_F$}

\noindent
The following theorem provides the estimate on the rate of convergence of the numerical method described in Subsection \ref{description}.

\begin{thm}\label{rec_error}
Let $\mu$ be a solution to \eqref{MainEq} with coefficients $(b,c,\beta,f)$ and initial data $\mu_o$. Let $\mu_{k}$ be a numerical solution at time $t_k = k\Delta t$
obtained by solving \eqref{subproblem1} - \eqref{subproblem2} with coefficients $(b,c,\beta,\bar f^{\epsilon})$ and initial data $\mu_o^{\delta}$, where $\bar f^{\epsilon}$ is the modified Finite Range Approximation of $f$ defined by \eqref{F_BARRR} and $\mu_o^{\delta}$ is a sum of Dirac deltas. Then, there exists a constant $C$, which depends on $(b,c,\beta,f)$, $\mu_o$, and $T$, such that
\begin{equation}\label{ppp}
\rho_F\left( \mu_{k},\mu(t_k) \right) \leq C\left( \Delta t +
(\Delta t)^{\alpha} + \mathcal{O}(\Delta t) \epsilon + \epsilon + \mathcal I(\mu_o) \right),
\end{equation}
where $\Delta t$ is the length of a time step, $\epsilon$ is the error of the modified Finite Range Approximation \eqref{F_BARRR}, that is $\norm{f - \bar f^{\epsilon}}_{\L\infty} < \epsilon$, and $\mathcal I(\mu_o) := \rho_F(\mu_o, \mu_o^{\delta})$ is the error of the initial data approximation.
\end{thm}

\begin{remark}\label{rmerror}
The error estimate \eqref{ppp} accounts for different error
sources. More specifically, the error of the order
$\mathcal{O}(\Delta t)$ is a consequence of the splitting
algorithm. The term of order $\mathcal{O}((\Delta t)^{\alpha})$
follows from the fact that we solve ODEs with
parameter functions independent of time, while $b,c$ and $\eta$
are in fact H\"older continuous with exponent $\alpha$ with respect to time.
The error of the initial data approximation $\mathcal I(\mu_o)$ is inversly proportional to the number of Dirac deltas approximating $\mu_o$. According to Lemma \ref{init}, $\mathcal I(\mu_o)$ can be arbitrarily small.
\end{remark}
\begin{proofof}{ Theorem \ref{rec_error}}
The proof is divided into several steps. For simplicity, in all
estimates below, we will use a generic constant $C$, without
specifying its exact form that may change from line to line.
\smallskip
\\
\textbf{Step 1: The auxiliary scheme.} \quad Let us define the
auxiliary semi-continuous scheme, which consists in solving
subsequently the problems
\begin{equation}\label{exact1}
\left\{
\begin{array}{rcl}
\frac{\partial}{\partial t} \mu + \frac{\partial}{\partial
x}(\bar b_k(x)\mu) &=& 0,\quad \mathrm{on}\;\; [t_k, t_{k+1}] \times \Rp
\\[1mm]
\mu(t_k) &=& \mu_k,
\end{array}
\right.
\end{equation}
and
\begin{equation}\label{exact2}
\left\{
\begin{array}{rcl}
\frac{\partial}{\partial t} \mu +\bar{\bar c}_k(x)\mu &=& \int_{\Rp} 
\bar{\bar\eta}_k(y)
d \mu_t(y),\quad \mathrm{on}\;\; [t_k, t_{k+1}] \times \Rp
\\[1mm]
\mu(t_k) &=& \bar\mu_k,
\end{array}
\right.
\end{equation}
where $\mu_k \in \mathcal M^+(\Rp)$, $\bar \mu_k$ is
the solution to \eqref{exact1} at time $t_{k+1}$ and $\bar b_k$,
$\bar{\bar c}_k$, and $\bar{\bar \eta}_k$ are defined as
\begin{eqnarray}
\label{bar_b_freeze}
 \bar b_k(x) &=& b \left(t_k,  \mu_{k}\right)(x),
\\[2mm]
\label{eta_zla}
\bar {\bar c}_k(x) &=& c \left(t_k,  \bar \mu_{k}\right)(x),
\quad
\bar{\bar \eta}_k(y) = \sum_{p=1}^{r}\beta_p(t_k,  \bar \mu_{k})(y) \; \delta_{x =  \bar f^{\epsilon}_p(y)}.
\end{eqnarray}
A solution to the second equation at time $t_{k+1}$ is denoted by $ \mu_{k+1}$. Denote by $\nu_{k+1}$ a solution to \eqref{exact2} with $\bar{\bar \eta}_k$ defined as $\sum_{p=1}^{r}\beta_p(t_k,  \bar \mu_{k})(y) \; \delta_{x =  f_p(y)}$.
\smallskip
\\
\textbf{Step 2: Error of the Finite Range Approximation.}\quad
According to \eqref{better_estimate_eq}, it holds that
    \begin{eqnarray}\label{osza}
    \rho_F\left(\mu_{k+1},\nu_{k+1}\right)
\leq
C\Delta t \e^{C \Delta t} \norm{f - \bar f^{\epsilon}}_{\L\infty}
\leq
     C\Delta t \e^{C \Delta t} \epsilon \leq \tilde C \epsilon \Delta t,
    \end{eqnarray}
    where $\tilde C$ is such that $C\e^{Ch} \leq \tilde C$ for all $h \in [0,T]$.
\smallskip
\\
\textbf{Step 3: Error of splitting.} \quad Let $\nu(t)$ be a
solution to \eqref{MainEq} on a time interval $[t_{k}, t_{k+1}]$ with
initial datum $\mu_{k}$ and parameter functions $\bar b_k$,
$\bar c_k$, $\bar \eta_k$, where $\bar b_k$ is defined by
\eqref{bar_b_freeze},
\begin{eqnarray}
\label{bar_c_freeze}
{{\bar c}}_k(x) &=& c \left(t_k,   \mu_{k}\right)(x),
\\
\label{bar_eta_freeze}
\bar{ {\eta}}_k(y) &=& \sum_{p=1}^{r}\bar{\beta}_{p}(t_k,\mu_k)(y) \; \delta_{x = f_p(y)}
= : \sum_{p=1}^{r}\bar{\beta}_{p,k}(y) \; \delta_{x = f_p(y)}.
\end{eqnarray}
According to \cite[Proposition 2.7]{ColomboGuerra2009} and
\cite[Proposition 2.7]{CCGU}, the distance between $ \nu_{k+1}$ and $\nu(t_{k+1})$, that is, the error coming from the
application of the splitting algorithm can be estimated as
\begin{equation}\label{est_split}
\rho_F( \nu_{k+1},\nu(t_{k+1})) \leq C (\Delta t)^2,
\end{equation}
where $C$ depends on the $\norm{\cdot}_{\W{1}{\infty}}$ norm of $b,c,\beta$ and the Lipschitz constant $\mathbf{Lip}(f_p)$, $p=1,\dots,r$.
\\[3mm]
\noindent
To estimate a distance between $\nu(t_{k+1})$ and $\mu(t_{k+1})$
consider $\zeta(t)$, which is a solution to \eqref{MainEq} on a time
interval $[t_k, t_{k+1}]$ with initial data $\mu(t_k)$ and
coefficients $\bar b_k$, $\bar c_k$, $\bar \eta_k$. By triangle
inequality
$$
\rho_F(\nu(t_{k+1}), \mu(t_{k+1})) \leq \rho_F(\nu(t_{k+1}), \zeta(t_{k+1})) + \rho_F(\zeta(t_{k+1}), \mu(t_{k+1})).
$$
The first term of the inequality above is a distance between
solutions to \eqref{MainEq} with different initial data, that is,
$ \mu_{k}$ and $\mu(t_k)$ respectively. The second term is
equal to a distance between solutions to \eqref{MainEq} with
coefficients $(\bar b_k, \bar c_k,\bar \eta_k)$ defined by
\eqref{bar_b_freeze}, \eqref{bar_c_freeze},
\eqref{bar_eta_freeze} and $(b(t,\mu(t)), c(t,\mu(t)),
\eta(t,\mu(t)))$, respectively. By the continuity of solutions to \eqref{MainEq}
with respect to the initial datum and coefficients in Theorem
\ref{better_estimate}, we obtain
\begin{eqnarray}\label{main1}
\rho_F(\nu(t_{k+1}), \zeta(t_{k+1}))
\leq
\e^{ C \Delta t} \rho_F( \mu_{k}, \mu(t_{k})),
\end{eqnarray}
and
\begin{align}\label{main2}
\rho_F(\zeta(t_{k+1}), & \mu(t_{k+1}))
 \leq C \Delta t  \e^{ C\Delta
t} \left( \norm{(\bar b_k - b,\bar c_k - c)}_{\overline{\BC}} + \sum_{p=1}^{r} \norm{\bar \beta_{p,k} -
{\beta_p}}_{\overline{\BC}} \right),
\end{align}
where
$\nonumber \norm{\bar f_k - f}_{\overline{\BC}} = \sup_{t \in
[t_k, t_{k+1}]}\norm{ \bar f_k - f(t,\mu(t))}_{\L\infty}$, and $f \in \{b, c, \beta_p\}$.
By the assumption \eqref{eqAssumptions} and definition \eqref{bar_b_freeze} of $\bar b_k$
\begin{eqnarray}
\nonumber
\norm{\bar b_k - b(t,\mu(t))}_{\L\infty}
&\leq&
\norm{ b(t_{k}, \mu_{k}) - b(t_{k},\mu(t))}_{\L\infty} +
\norm{b(t_{k},\mu(t)) - b(t,\mu(t))}_{\L\infty}
\\
&\leq&
 \mathbf{Lip}(b(t_k, \cdot))\; \rho_F( \mu_{k}, \mu(t)) +
 \f{b} \modulo{t - t_k}^{\alpha}.
\label{333}
\end{eqnarray}
Using the Lipschitz continuity of the solution $\mu(t)$ (Theorem \ref{better_estimate}) yields
$$
\rho_F( \mu_{k}, \mu(t))\leq
 \rho_F( \mu_{k}, \mu(t_k)) + \rho_F(\mu(t_k), \mu(t))
\leq
 \rho_F( \mu_{k}, \mu(t_k))  + C \Delta t \e^{ C \Delta t}.
$$
Substituting the latter expression into \eqref{333} yields
$$
\norm{\bar b_k - b(t,\mu(t))}_{\L\infty} \leq
\mathbf{Lip}(b)\left(
 \rho_F( \mu_{k}, \mu(t_k))  + C \Delta t \e^{ C \Delta t}
\right)
+\f{b}(\Delta t)^{\alpha},
$$
where $\mathbf{Lip}(b) = \sup_{t\in[0,T]}  \mathbf{Lip}(b(t, \cdot))$.
Bounds for $\norm{\bar c_{k} -
{c}}_{\overline{\BC}}$ and $\norm{\bar \beta_{p,k} -
{\beta_p}}_{\overline{\BC}}$ can be proved analogously.
From the assumptions it holds that
$$
\norm{(b,c,\beta)} : = \mathbf{Lip}(b) + \mathbf{Lip}(c) + \f{(b,c)} + \sum_{p=1}^{r}\left(
\mathbf{Lip}(\beta_p) + \f{\beta_p}
\right)< +\infty,
$$
and as a consequence, we obtain
\begin{align*}
\norm{\bar b_k - b}_{\overline{\BC}} + \norm{\bar c_k -
c}_{\overline{\BC}} & +\sum_{p=1}^{r} \norm{\bar \beta_{p,k} -
{\beta_p}}_{\overline{\BC}} 
\\
&\leq
\norm{(b,c,\beta)} \left( \rho_F(
\mu_{k}, \mu(t_k)) + C\Delta t \e^{ C  \Delta t}  + (\Delta
t)^{\alpha}\right).
\end{align*}
Using this inequality in \eqref{main2} yields
\begin{eqnarray*}
\rho_F(\zeta(t_{k+1}),  \mu(t_{k+1})) &\leq& C \Delta t  \e^{
C\Delta t} \Big(
 \rho_F( \mu_{k}, \mu(t_k))
+
 \Delta t
+
 (\Delta t)^{\alpha} \Big)
\\
&\leq&
C \Delta t  \e^{ C\Delta t}
 \rho_F( \mu_{k}, \mu(t_k))
+ C   \e^{ C T} (\Delta t)^2 + C \e^{ C T} (\Delta
t)^{1+\alpha}.
\end{eqnarray*}
Combining the inequality above with \eqref{main1} and redefining
$C$ leads to
\begin{eqnarray}\label{z2}
\rho_F(\nu(t_{k+1}), \mu(t_{k+1})) &\leq& \e^{ C \Delta t}(1 +
C\Delta t) \rho_F( \mu_{k}, \mu(t_{k})) + C (\Delta t)^2 + C
(\Delta t)^{1+\alpha} 
\nonumber
\\
 &\leq& \e^{ 2C \Delta t} \rho_F( \mu_{k},
\mu(t_{k})) + C (\Delta t)^2 + C (\Delta t)^{1+\alpha}.
\end{eqnarray}
Finally, putting together \eqref{osza}, \eqref{est_split}, and \eqref{z2} we
conclude that
\begin{equation}\label{znew}
\rho_F\left( \mu_{k+1}, \mu(t_{k+1})\right) \leq \e^{ 2C \Delta t}
\rho_F( \mu_{k}, \mu(t_{k})) + C (\Delta t)^2 + C (\Delta
t)^{1+\alpha} + C\epsilon\Delta t.
\end{equation}
\noindent
\textbf{Step 4: Adding the errors.} Application
of the discrete Gronwall's inequality to \eqref{znew} yields
\begin{eqnarray*}
\rho_F( \mu_{k}, \mu(t_{k})) 
&\leq&
 \e^{Ck\Delta t}  \rho_F( \mu_o^{\delta}, \mu_o)
+ C \frac{\e^{C k\Delta t}-1}{\e^{C\Delta t} - 1}  \left( (\Delta
t)^2 + (\Delta t)^{1+\alpha} + \epsilon\Delta t \right).
\end{eqnarray*}
There exists a constant $C^*$ depending only on $T$ such that $\e^{Ck \Delta
t}-1 <  C^*k \Delta t$, for each $k \Delta t \in [0,T]$.
Therefore, we deduce
$$
\frac{\e^{C k\Delta t}-1}{\e^{C\Delta t} - 1} \leq
\frac{C^*k\Delta t}{C \Delta t} = \frac{C^*}{C}k,
$$
and thus,
\begin{equation*}
\rho_F( \mu_{k}, \mu(t_{k})) 
\leq
 \e^{Ck\Delta t}  \mathcal I(\mu_o)
+ C k\Delta t \left( \Delta t + (\Delta t)^{\alpha} +  \epsilon\right).
\end{equation*}
Since
$k \Delta t = t_k$, the assertion is proved.
\end{proofof}
\noindent
In the following lemma we show that any measure $\nu \in \Mp$ can be approximated in $\rho_F$ with an arbitrarily small error by a sum of Dirac deltas.
\begin{lemm}\label{init}
Let $\nu \in \Mp$ be such that $M_{\nu} = \int_{\Rp} \d \nu \neq 0$. Then, for each $\delta > 0$ there exists $K \in \mathbbm{N}$ and a measure $\tilde \nu = \sum_{i=1}^{K} m_i \delta_{x_i}$, such that
\begin{eqnarray}
\rho_F(\nu,\tilde \nu)\leq \delta.
\end{eqnarray}
\end{lemm}

\begin{proofof}{Lemma \ref{init}}
A measure $\nu \in \Mp$ is tight. Therefore, for each $\epsilon > 0$ there exists $K_{\epsilon} > 0 $ such that $\nu(\Rp \backslash [0,K_{\epsilon}]) \leq \epsilon$. Define $\nu^{\epsilon}$ as a restriction of $\nu$ to $[0,K_{\epsilon}]$. Let $\phi \in \W{1}{\infty}(\Rp)$. Then,
\begin{eqnarray*}
\int_{\Rp}\phi (x)\d (\nu - \nu^{\epsilon})(x) &=&
\int_{[0,K_{\epsilon}]}\phi (x)\d (\nu - \nu^{\epsilon})(x)  + \int_{(K_{\epsilon}, + \infty)}\phi (x)\d (\nu - \nu^{\epsilon})(x) 
\\
&=&
 \int_{(K_{\epsilon}, + \infty)}\phi (x)\d \nu (x) \leq \norm{\phi}_{\L\infty} \epsilon.
\end{eqnarray*} 
Taking supremum over all $\phi \in \W{1}{\infty}(\Rp)$ such that $\norm{\phi}_{\L\infty} \leq 1$ yields
\begin{equation}\label{q1}
\rho_F(\nu, \nu^{\epsilon}) \leq \epsilon.
\end{equation}
Let $M_{\nu}^{\epsilon} = \int_{\Rp} \d \nu^{\epsilon} = \int_{[0,K_{\epsilon}]} \d \nu$. Then, according to \cite[Lemma 2.1 and (2.17)]{CGU}, the error of the fixed-location approximation of $\nu^{\epsilon}$ by a measure $\tilde \nu$ consisting of $Q$ Dirac deltas is equal to
\begin{equation}\label{q2}
\rho_F(\nu_{\epsilon}, \tilde \nu) \leq M_{v}^{\epsilon} \frac{K_{\epsilon}}{2Q}.
\end{equation}
Taking $\displaystyle \epsilon = \frac{\delta}{2}$ in \eqref{q1} and $ \displaystyle Q \geq \frac{4}{M_{\nu}^{\epsilon}K_{\epsilon} \delta}$ in \eqref{q2} finishes the proof.
\end{proofof}

\section{Numerical Results}\label{num}

The aim of this section is to perform the particle-based method described in Subsection \ref{sec_opis}, and present results of the simulations for the symmetric cell division model. In this test case, a birth process occurs due to a division of a mother cell into two cells of equal sizes. Moreover, we assume the existence of both, minimal and maximal cell reproduction sizes, $x_o>0$ and $x_{max} > x_o$ respectively. Namely, a single cell is not able to divide, unless it reaches size $x_o$, and it divides before reaching its maximal size $x_{max}$ with probability equal to $1$. From this assumption it follows the existence of the minimal cell size, which is equal to $\frac{1}{2}x_o$. Following \cite{abia}, we set the coefficients as below
\begin{eqnarray*}
&&x_0=\frac 1 4, \quad x_{max}=1,
\\
&&b(x)=0.1(1-x), \quad c(x)= \beta(x), \quad \eta(t,\mu)(y)=2\beta(y)\delta_{x=\frac 1 2 y},
\end{eqnarray*}
where
$$\beta(y)=\begin{cases}0, &\mbox{ for } y\in(-\infty,x_0)\cup(x_{max},+\infty),\\ \frac{b(y)g(y)}{1-\int_{x_0}^y g(x)dx}, &\mbox{ for } y\in[x_0,x_{max}],  \end{cases}$$
and 
$$g(y)=\begin{cases} \frac{160}{117}(-\frac 2 3 +\frac 8 3 y)^3, &\!\!\!\!\mbox{ for } y\in[x_0,\frac{x_0+1}{2}],\\
\frac{640}{117}\naw{-1+2y+\frac{16}{3}(y-\frac 5 8)^2}+ \frac{5120}{9}(y-\frac 5 8)^3(\frac 8 3 y - \frac{11}{3}), &\!\!\!\!\mbox{ for } y\in(\frac{x_0+1}{2},x_{max}].  \end{cases}$$
We consider the initial data given by $$\mu_0(x)=(1-x)(x-\frac 1 2 x_0)^3.$$
In our case, the exact solution is unknown. Therefore, we calculate a reference solution $\mu_{ref}$ which is a numerical solution with the following (small) parameters
$$\epsilon~=~3.90625\cdot 10^{-05}, \quad \mathrm{and}\quad \Delta t=\Delta x=0.00078125.$$
We recall that $\epsilon$ denotes the accuracy of the Finite Range Approximation, see Lemma~\ref{ccc}.  According to our leading assumption \eqref{form_of_eta}, $\eta$ has the following form
$$\eta(t,\mu)(y) = \sum_{p=1}^{r}\beta_p(t,\mu)(y) \delta_{x = f_p(y)}.$$
In case of the symmetric cell division we set $p=1$ and $f_p(y) = \frac{1}{2}y$. $\Delta t$ denotes the length of a time step, and $(x_{max} - \frac{1}{2}x_o)/\Delta x$ is a number of Dirac deltas approximating the initial data $\mu_o$. This number is inversly proportional to the error of approximation of the initial data, see Lemma \ref{rmerror}.
The error is defined by the following formula
$$\Err(T, \Delta t, \epsilon) := \rho(\mu_{ref},\mu_{{k}}),$$ where $k$ is such that $k\Delta t = T$, and $T$ is the final time. Function $\rho$ is defined as
$$\rho(\mu,\nu)=\min\{M_\mu,M_\nu\}W_1\naw{\frac{\mu}{M_\mu},\frac{\nu}{M_\nu}}+|M_\mu-M_\nu|,$$
where $M_\mu=\int_{\Rp} d\mu$, and $W_1$ is the 1-Wasserstien metric on the space of probability measures, which in the one dimensional case can be obtained by the formula
$$W_1(\mu_1,\mu_2)=\int_{\Rp}|F_{\mu_1}(x)-F_{\mu_2}(x)|dx,$$
where $F_{\mu_i}$ denotes the distribution function of the measure $\mu_i$. The function $\rho$ is equivalent to the flat metric $\rho_F$ in the sense that there exists a constant $C$ such that
$$
C \rho(\mu_1, \mu_2) \leq \rho_F(\mu_1,\mu_2) \leq \rho(\mu_1,\mu_2),
$$
see \cite[Lemma 2.1]{CGU} for details. The order of the method is given by
$$q:=\lim_{\Delta t \to 0} \frac{\log\naw{\frac{\Err(T,2\Delta t,\epsilon)}{\Err(T,\Delta t,\epsilon)}}}{\log 2}.$$
In Tables \ref{tabela1} - \ref{tabela3} we present results of numerical simulations with different values of $\Delta t$ and $\epsilon$. Parameter $\Delta x$ is always equal to the corresponding value of $\Delta t$. 
\\
~\\
\noindent
Conclusions of our numerical tests are the following. For a fixed value of $\epsilon$ the order of convergence tends to zero as $\Delta t = \Delta x << \epsilon$, see Table \ref{tabela1}. It is intuitively clear, since $\epsilon$ denotes the accuracy of the Finite Range Approximation of $\eta$. Therefore, decreasing the time step $\Delta t$ and the parameter $\Delta x$ does not decrease the error, since $\eta$ is not sufficiently accurately approximated. When $\epsilon$ is relatively small ($10^{-3}$ or $10^{-4}$) the order of convergence tends to $1$, see Table \ref{tabela2} and Table \ref{tabela3}. Summing up, the optional approach is to use $\epsilon$ proportional to $\Delta t$ and $\Delta x$, which is consistent with the theoretical error estimate from Theorem \ref{rec_error}.
\\
~\\
\begin{table}[ht]
\begin{center}
\begin{tabular}{ |c|c|c| }
  \hline
  $\Delta t$ & $\Err(1,\Delta t,10^{-2})$ & $q$ \\ \hline
 $1.0000\cdot 10^{-1}$ & $1.039119349267744\cdot 10^{-3}$ & -\\
 $5.0000\cdot 10^{-2}$ & $8.823037689904696\cdot 10^{-4}$ & 0.2360140132647169\\
 $2.5000\cdot 10^{-2}$ & $1.559106188788983\cdot 10^{-4}$ & 2.5005562562550110\\
 $1.2500\cdot 10^{-2}$ & $8.584324175849732\cdot 10^{-5}$ & 0.8609427278124292\\
 $6.2500\cdot 10^{-3}$ & $5.219135515056872\cdot 10^{-5}$ & 0.7178936965493665\\
 $3.1250\cdot 10^{-3}$ & $3.770534034670716\cdot 10^{-5}$ & 0.4690419898814134\\
 $1.5625\cdot 10^{-3}$ & $3.055577759898044\cdot 10^{-5}$ & 0.3033236762664360\\
 $7.8125\cdot 10^{-4}$ & $2.705291385881350\cdot 10^{-5}$ & 0.1756612013876553\\
  \hline
\end{tabular}
  \caption{The error of the finite range method for $\epsilon=10^{-2}$ on the interval $[0,1]$}\label{tabela1}
\end{center}
\end{table}

\begin{table}[ht]
\begin{center}
\begin{tabular}{ |c|c|c| }
  \hline
  $\Delta t$ & $\Err(1,\Delta t,10^{-3})$ & $q$ \\ \hline
 $1.0000\cdot 10^{-1}$ & $1.023703008152967\cdot 10^{-3}$ & -\\
 $5.0000\cdot 10^{-2}$ & $8.618079480449228\cdot 10^{-4}$ & 0.2483589193046122\\
 $2.5000\cdot 10^{-2}$ & $1.446641895315060\cdot 10^{-4}$ & 2.5746585654081700\\
 $1.2500\cdot 10^{-2}$ & $7.307539446483757\cdot 10^{-5}$ & 0.9852502214828665\\
 $6.2500\cdot 10^{-3}$ & $3.679026124907980\cdot 10^{-5}$ & 0.9900617908153437\\
 $3.1250\cdot 10^{-3}$ & $1.870962335325300\cdot 10^{-5}$ & 0.9755434052597932\\
 $1.5625\cdot 10^{-3}$ & $9.718941155567934\cdot 10^{-6}$ & 0.9449094647632132\\
 $7.8125\cdot 10^{-4}$ & $5.916312147904012\cdot 10^{-6}$ & 0.7161009740113241\\
  \hline
\end{tabular}
  \caption{The error of the finite range method for $\epsilon=10^{-3}$ on the interval $[0,1]$}\label{tabela2}
\end{center}
\end{table}

\begin{table}[ht]
\begin{center}
\begin{tabular}{ |c|c|c| }
  \hline
  $\Delta t$ & $\Err(1,\Delta t,10^{-4})$ & $q$ \\ \hline
 $1.0000\cdot 10^{-1}$ & $1.021598713535499\cdot 10^{-3}$ & - \\
 $5.0000\cdot 10^{-2}$ & $8.599090072320180\cdot 10^{-4}$ & 0.2485727017000523\\
 $2.5000\cdot 10^{-2}$ & $1.441428362080116\cdot 10^{-4}$ & 2.5766848688633300\\
 $1.2500\cdot 10^{-2}$ & $7.245219715669196\cdot 10^{-5}$ & 0.9923977915838431\\
 $6.2500\cdot 10^{-3}$ & $3.615804743340988\cdot 10^{-5}$ & 1.0027126679987400\\
 $3.1250\cdot 10^{-3}$ & $1.787395766731422\cdot 10^{-5}$ & 1.0164576605004100\\
 $1.5625\cdot 10^{-3}$ & $8.715707767967237\cdot 10^{-6}$ & 1.0361693827983020\\
 $7.8125\cdot 10^{-4}$ & $4.150900779265493\cdot 10^{-6}$ & 1.0701933776302790\\
  \hline
\end{tabular}
  \caption{The error of the finite range method for $\epsilon=10^{-4}$ on the interval $[0,1]$}\label{tabela3}
\end{center}
\end{table}


\section{Appendix}

\begin{proofof}{Lemma \ref{Lema1}}
Let $\Psi \in (\W{1}{\infty}(\Rp))^*$, $\phi \in \W{1}{\infty}(\Rp)$ be such that $\norm{\Psi}_{(\W{1}{\infty})^*} \leq 1$ and $\norm{\phi}_{\W{1}{\infty}} \leq 1$. Then,
\begin{equation}\label{szacowanie}
\begin{aligned}
\norm{T_t - I}_{L\naw{(\W{1}{\infty})^*}} 
& =\sup_{\Psi}
 \norm{T_t\Psi - \Psi}_{(\W{1}{\infty})^*}
\\
& = \sup_{\Psi}\; \sup_{\phi} \abs{T_t\Psi(\phi) - \Psi(\phi)}
\geq \abs{T_t\Psi_\lambda(\bar \phi) - \Psi_\lambda(\bar \phi)},
\end{aligned}
\end{equation}
where $\Psi_{\lambda} \in (\W{1}{\infty}(\Rp))^*$ and $\bar \phi\in \W{1}{\infty}(\Rp)$ are arbitrary. In what follows we construct $\Psi_{\lambda}$ and $\bar\phi$, such that $\abs{T_t\Psi_\lambda(\bar \phi) - \Psi_\lambda(\bar \phi)} = 1$.
\\[2mm]
\noindent
Let $\delta_1$ be a Dirac measure located in $x = 1$. In fact, $\delta_1$ is the functional on the space of bounded continuous functions $(\mathbf{BC}(\Rp), \norm{\cdot}_{\infty})$ such that $\delta_1(f) = f(1)$. This functional can be extended to the linear space
$$X=\{f\colon\Rp\to\RR\   :   f \in \mathbf{B}(\Rp) \mbox{ and $f$ has one-sided limits in } x=1\} \subset \L{\infty}(\Rp),$$
where $\mathbf{B}(\Rp)$ denotes the space of measurable and bounded functions on $\Rp$. This extension is obtained by putting
$${\delta}_1^\lambda (f) = \lim_{r\to 0} \left(\lambda f(1-r)+(1-\lambda)f(1+r)\right),
$$
where $0 \leq \lambda \leq 1$.
It holds that
$$\abs{{\delta}_1^\lambda(f)} \leq \norm{f}_{\L{\infty}}.$$
The Hahn-Banach extension theorem guaratees existence of an extension of ${\delta}_1^\lambda$ to $\L{\infty}(\Rp)$, which is further denoted as $\bar{\delta}_1^\lambda$. 
Therefore, setting
\begin{equation}\label{Del}
\Psi_\lambda(\phi) : = \bar{\delta}_1^\lambda(\phi'),
\end{equation}
where $\phi'$ denotes the derivative of $\phi$,
yields a family of bounded linear functionals $\Psi_\lambda$ on $\W{1}{\infty}$ with $\norm{\Psi_\lambda}_{(\W{1}{\infty})^*}=1$. 
Let
\begin{equation*}
f(x) = \begin{cases} 1, & x\in[0,1]\\ 0, & x>1  \end{cases}
\quad\quad
\mathrm{and}
\quad\quad
\bar \phi(x) = \int_{0}^{x}f(t) dt.
\end{equation*}
Clearly, $f\in \L{\infty}(\Rp)$, $\bar\phi\in\W{1}{\infty}(\Rp)$ and $\norm{\bar \phi}_{\W{1}{\infty}} = 1$.
Inserting $\Psi_{\lambda}$ defined by \eqref{Del} and $\bar \phi$ defined as above into \eqref{szacowanie} yields
$$\begin{aligned}
 \Big| (T_t \Psi_\lambda & - \Psi_\lambda)(\bar\phi) \Big| =
 \abs{\naw{T_t \Psi_\lambda - \Psi_\lambda}(\bar\phi)} = 
\abs{\Psi_\lambda(T_{-t}\bar\phi - \bar\phi)} 
\\
& = \abs{\bar{\delta}_1^\lambda((T_{-t}\bar\phi - \bar\phi)')}  = \abs{\bar{\delta}_1^\lambda((\bar\phi(\cdot + t) - \bar\phi(\cdot))')} = \abs{{\delta}_1^\lambda(f(\cdot + t) - f(\cdot))} =
\\[2mm]
& = \lim_{r\to 0^+} \abs{ \lambda\naw{f(1-r+t)- f(1-r)}+(1-\lambda)\naw{f(1+r+t)-f(1+r)}}  = \\
& = \lambda  \lim_{r\to 0^+}  \abs{f(1-r+t)- f(1-r)} = \lambda,
\end{aligned}$$
which implies
\begin{equation*}
\norm{T_t - I}_{L\naw{(\W{1}{\infty})^*}} \geq  \abs{T_t\Psi_\lambda(\bar\phi) - \Psi_\lambda(\bar\phi)} = \lambda.
\end{equation*}
Taking supremum over $\lambda \in [0,1]$ finishes the proof.
\end{proofof}

\section*{Acknowledgements}

{\small PG and AU were supported by the International PhD Projects Programme of Foundation for Polish Science operated within the Innovative Economy Operational Programme 2007-2013 (PhD Programme: Mathematical Methods in Natural Sciences). PG was also supported by the grant IdP2011/000661.}

\end{document}